\documentclass[preprint]{amsart}
\usepackage{amsrefs}

\newcommand{\R}{\mathbf R}

\newcommand{\grad}{\operatorname {grad}}

\newtheorem{theorem}{Theorem}
\newtheorem{lemma}{Lemma}
\newtheorem{corollary}{Corollary}

\newcommand{\supp}{\operatorname {supp}}
\newcommand{\norm}[1]{\lVert#1\rVert}

\newcommand{\beq}[1]{\begin{equation}\label{#1}}
\newcommand{\eeq}{\end{equation}}

\newcommand{\beqan}{\begin{eqnarray}}
\newcommand{\eeqan}{\end{eqnarray}}

\title{The Hohenberg-Kohn Theorem for Schrodinger Semigroups}

\author{Omar Hijab}
\address{Department of Mathematics, Temple University, Philadelphia, PA 19122}
\email{hijab@temple.edu}


\date{\today}

\begin{document}

\maketitle

\begin{abstract}
At the basis of much of computational chemistry is density functional theory,
as initiated by the Hohenberg-Kohn theorem. The theorem states that, when nuclei are fixed,
electronic systems are determined by $1$-electron densities.  We recast and derive this result
within the context of the principal eigenvalue of Schrodinger semigroups.
\end{abstract}

\section{Introduction}

In quantum mechanics, the probability distribution of the ground state of an $N$-electron system\footnote{An atom, molecule, or solid where nuclei are fixed.}
is a permutation-symmetric probability measure $\mu$ on $\R^{3N}$, and its $1$-electron marginal 
is the probability measure $\rho$ on $\R^3$ given by

$$\int_{\R^3}f\,d\rho = \int_{\R^{3N}} f(x_1)\,d\mu(x_1,\dots,x_N).$$
The potential acting on the electrons is a sum $V_0+V$ of potentials, where $V_0$ is the repulsive 
Coulomb potential between electrons, and $V$ is the attractive nuclear or external potential\footnote{The $1/N$ normalization is not standard.}

\beq{eq:ext}
V(x_1,\dots,x_N) = \frac{v(x_1)+\dots+v(x_N)}N,
\eeq
for some function $v$ on $\R^3$. The system is specified by the external potential $v$, as $V_0$ is the same for all $N$-electron systems.

Then the electronic ground state energy is given by
\beq{eq:atom}
E(V_0+V) = \inf_\psi \int_{\R^{3N}}\left(|\grad\psi|^2 + V_0\psi^2+V\psi^2\right)\,dx_1\dots dx_N,
\eeq
where the infimum is over all real $\psi$ satisfying $\int\psi^2dx_1\dots dx_N=1$, and the
distribution corresponding to the ground state $\psi$ is $d\mu = \psi^2\,dx_1\dots dx_N$.

The Hohenberg-Kohn theorem \cite{HK} states that the external potential $v$ --- and thus the electronic system --- is determined by the marginal
$\rho$: If $\mu_1$, $\mu_2$ are distributions of ground states $\psi_1$, $\psi_2$ corresponding to external potentials
$v_1$, $v_2$, and their marginals agree, $\rho_1=\rho_2$, then $v_1-v_2$ is a constant. The thrust of the theorem  is to reduce the study of electronic systems from $3N$ variables down
to $3$ variables.

In this paper we generalize this result from the above electronic setting to the general (non-self-adjoint)
Markov semigroup setting. To help simplify matters, instead of $\R^3$, we take a compact metric space $X$
as our position space.

Let $X$ be a compact metric space and let $P_t$, $t\ge0$, be a Markov semigroup on $C(X)$ with 
generator $L$ defined on its dense domain $\mathcal D\subset C(X)$. Examples of semigroups which satisfy 
all our assumptions below are

\begin{itemize}

\item $X$ is a compact manifold and $L$ is a nondegenerate elliptic second order differential operator 
with smooth coefficients, given by

$$Lf(x) = \sum a_{ij}(x)\frac{\partial^2f}{\partial x_i\partial x_j} + \sum b_i(x)\frac{\partial f}{\partial x_i}$$
in local coordinates. 

\item $X=\{1,\dots,d\}$ and $L$ is a $d\times d$ matrix with nonnegative off-diagonal entries whose
row-sums vanish and whose adjacency graph is connected.

\end{itemize}

Given $V$ in $C(X)$, let $P_t^V$, $t\ge0$, denote the Schrodinger semigroup on $C(X)$ generated by $L+V$. 
Then the principal eigenvalue

$$\lambda_V\equiv \lim_{t\uparrow\infty} \frac1t \log\norm{P_t^V}$$
exists and is given by the Donsker-Varadhan formula \cite{DV}

\beq{eq:DV}
\lambda_V = \sup_\mu \left(\int_X V\,d\mu - I(\mu)\right)
\eeq
where the supremum is over all probability measures $\mu$ on $X$ and 

$$I(\mu) = - \inf_{u\in\mathcal D^+}\int_X\frac{Lu}{u}\,d\mu.$$
Here the infimum is over all positive $u$ in $\mathcal D$. 
In the electronic case, (\ref{eq:DV}) reduces to (\ref{eq:atom}) and $\lambda_V = -E(-V)$.

Given $f\in C(X)$ and a probability measure $\mu$ on $X$, let $\mu(f)$ denote the integral of $f$ against $\mu$. 
Let $M(X)$ denote the space of probability measures on $X$, and let $V$ be in $C(X)$. 

An {\em equilibrium measure for $V$} is a $\mu\in M(X)$ achieving\footnote{The supremum is always achieved as $I$ is lower semicontinuous (Lemma \ref{lemma:lsc}).}
 the supremum  in (\ref{eq:DV}), $\lambda_V=\mu(V)-I(\mu)$.

A {\em ground measure for $V$} is a $\pi\in M(X)$ satisfying
\beq{eq:ground}
\int_Xe^{-\lambda_Vt}P^V_tf\,d\pi = \int_X f\,d\pi,\qquad t\ge0, f\in C(X).
\eeq

By positivity,
\beq{eq:posmeas}
P_t^Vf(x) = \int_X p^V(t,x,dy) f(y)
\eeq
for some family $(t,x)\mapsto p^V(t,x,\cdot)$ of bounded positive measures on $X$.
Thus $0\le P_t^Vf(x)\le +\infty$ is well-defined for $f$ nonnegative Borel on $X$.
Let $\mu$ be in $M(X)$. 

A {\em ground state for $V$ relative to $\mu$} is a nonnegative 
Borel function $\psi$ on $X$ satisfying $\psi>0$ a.e. $\mu$ and
$$e^{-\lambda_Vt}P^V_t\psi = \psi,\qquad a.e. \mu, t\ge0.$$ 

Thus a ground state $\psi$ plays the role of a right eigenvector for $L+V$, and a ground measure $\pi$ 
plays the role of a left eigenvector for $L+V$, both with eigenvalue $\lambda_V$. 

When $N=1$, the Hohenberg-Kohn theorem states that if $\mu$ is the distribution of a ground state $\psi$ 
corresponding to $V_1$ and to $V_2$, then $V_1-V_2$ is a constant. In the electronic case, $d\mu = \psi^2\,dx$ and this is an immediate consequence of the Schrodinger equations $L\psi+ V_i\psi = \lambda_{V_i}\psi$, $i=1,2$. 
In the general case, however, establishing this turns out to be the heart of the matter, as the correspondence 
between equilibrium measures $\mu$ and ground states $\psi$ is not as direct. The following sheds
light on the relation between $\mu$, $\psi$, and $\pi$.

\begin{theorem}
\label{theorem:equiv}
Let $\mu,\pi\in M(X)$ and let $V\in C(X)$. 
Suppose $\mu<<\pi$ and suppose $\psi=d\mu/d\pi$ satisfies $\log\psi\in L^1(\mu)$.  Then the following hold.
\begin{itemize}
\item If $\pi$ is a ground measure for $V$ and $\psi$ is a ground state for $V$ relative to $\mu$,  
then $\mu$ is an equilibrium measure for $V$.
\item If $\pi$ is a ground measure for $V$ and $\mu$ is an equilibrium measure for $V$,  
then $\psi$ is a ground state for $V$ relative to $\mu$.
\item If $\mu$ is an equilibrium measure for $V$ and $\psi$ is a ground state for $V$ relative to $\mu$,  
then $\pi$ is a ground measure for $V$.
\end{itemize}
\end{theorem}

In the electronic case, $L$ is self-adjoint relative to $dx_1\dots dx_N$,  so heuristically a right eigenvector is a left eigenvector, so  a ground state $\psi$ leads to a ground measure $d\pi=\psi\,dx_1\dots dx_N$ and to an equilibrium measure $d\mu = \psi\,d\pi = \psi^2\,dx_1\dots dx_N$.

Given $\psi$ nonnegative, let
\beq{eq:Ppsi}
P^{V,\psi}_tf=\frac{e^{-\lambda_Vt}P^V_t(f\psi)}{\psi}.
\eeq
Then $P^{V,\psi}_tf(x)$ is defined at a point $x$ if $P^V_t(|f|\psi)(x)<\infty$ and $\psi(x)>0$.

\begin{theorem}
\label{theorem:pf}
Fix $V\in C(X)$ and suppose
\beq{eq:supbound}
C\equiv \sup_{t\ge0} \left( e^{-\lambda_Vt}\norm{P^V_t} \right) < \infty,
\eeq
and let $\mu$ be an equilibrium measure for $V$. Then there is a ground state $\psi$ for $V$ relative 
to $\mu$ and a ground measure $\pi$ for $V$ such that 
\begin{itemize}
\item $\log\psi\in L^1(\mu)$,
\item $\mu<<\pi$ and $d\mu/d\pi = \psi$, and
\item $P^{V,\psi}_t$, $t\ge0$,  is a Markov semigroup on $L^1(\mu)$, and
$\mu$ is $P^{V,\psi}_t$, $t\ge0$,  invariant
$$\int_XP^{V,\psi}_tf\,d\mu = \int_X f\,d\mu,\qquad f\in L^1(\mu), t\ge0.$$
\end{itemize}
\end{theorem}

Note this existence result is not just a Perron-Frobenius result, as $\psi$ and $\pi$ are
determined subordinate to the given equilibrium measure $\mu$.

Now we list our assumptions on the Markov semigroup $P_t$, $t\ge0$.

We assume a strong uniformity condition 

\begin{itemize}

\item[(A)] There is a $T>0$ and an $\epsilon=\epsilon(T)>0$ such that $P_T|f|(x)\ge \epsilon P_T|f|(y)$ for
all $x,y\in X$ and $f\in C(X)$.

\end{itemize}

As we shall see, (A) implies (\ref{eq:supbound}). We also assume

\begin{itemize}

\item[(B)] There is a $T>0$ such that $f\ge0$ in $C(X)$ implies $P_Tf>0$ everywhere in $X$. 

\end{itemize}

A core for $P_t$, $t\ge0$, is a subspace 
$\mathcal D^\infty\subset\mathcal D$ whose closure in the graph norm $\norm{f}+\norm{Lf}$  equals $\mathcal D$. 
We assume

\begin{itemize}

\item[(C)] There is a core $\mathcal D^\infty$ that is closed under multiplication and division: If 
$f,g\in \mathcal D^\infty$ then $fg\in \mathcal D^\infty$, and if moreover $g>0$, then $f/g\in\mathcal D^\infty$.

\end{itemize}

The square-field operator is

$$\Gamma(g) = L(g^2) - 2gLg,\qquad g\in\mathcal D^\infty.$$

Let $p(t,x,dy) = p^0(t,x,dy)$. As we have

$$\Gamma(g)(x) = \lim_{t\to0} \frac1t \int_X p(t,x,dy)\left(g(y)-g(x)\right)^2,$$
it follows that $\Gamma(g)\ge0$ for $g\in\mathcal D^\infty$. Below in Lemma \ref{lemma:core}, we show\footnote{This reduces to the definition of $\Gamma(g)$
when $f=1$.} 
\beq{eq:gammaineq}
\max f\cdot \Gamma(g) \ge L(fg^2)-2gL(fg)+g^2Lf \ge \min f \cdot \Gamma(g)
\eeq
for $f,g\in\mathcal D^\infty$.
We assume the nondegeneracy condition

\begin{itemize}

\item[(D)] If $g\in\mathcal D^\infty$ and $\Gamma(g)\equiv0$, then $g$ is a constant.

\end{itemize}

Let $B(X)$ denote the bounded Borel functions on $X$. 
We say a potential $V$ is {\em smooth} if $P^V_t$ maps  $B(X)$ into $\mathcal D^\infty$ for $t>0$.  
This depends on both $L$ and $V$. 

For the examples above, (A) and (B) are valid, and (C) and (D) are valid if we take $\mathcal D^\infty=C^\infty(X)$,  and $V$ is smooth in the above sense if $V$ is in $C^\infty(X)$
(for the second example, $C^\infty(X)=C(X)=B(X)$ equals all functions on $X$).

\begin{theorem}
\label{thm:N=1}
Assume (A), (B), (C), (D) and let $V_1,V_2$ be smooth potentials.
If $\mu$ is an equilibrium measure for $V_1$ and for $V_2$, then $V_1-V_2$ is a constant.
\end{theorem}

This result should hold more broadly, in which case one should obtain $V_1-V_2$ is a constant 
on the support of $\mu$. This restriction is natural because one cannot expect to determine the potential in regions outside the electron cloud. The more general  result is easily verified when $L\equiv0$ for any $V_1,V_2\in C(X)$, so nondegeneracy should not play a role in a broader formulation. A discrete time version of Theorem \ref{thm:N=1} in the case $X=\{1,\dots,d\}$ is in \cite{F}.

Note that $\mu$ is an equilibrium measure for $V$ iff $V$ is a subdifferential of $I$ at $\mu$, i.e. iff
$$I(\nu) \ge I(\mu) + \nu(V) - \mu(V),\qquad \nu\in M(X).$$
Subdifferentials at a given $\mu$ need not exist.  When subdifferentials do exist, Theorem \ref{thm:N=1} provides conditions under which uniqueness holds at the given $\mu$, up to a constant.

Next we look at Markov semigroups on $C(X^N)$.

Let $N\ge1$ and $X^N$ be the $N$-fold product of $X$. Let $P_t$, $t\ge0$, be a Markov semigroup on $C(X^N)$, 
representing the motion of $N$  particles, and let $L$ be its generator. Let $P^i_t$, $t\ge0$, $1\le i\le N$, 
be Markov semigroups on $C(X)$. When $P_t$, $t\ge0$,
is the product of  $P^i_t$, $t\ge0$, $1\le i\le N$, with the $i$-th semigroup acting on the
$i$-th component in $C(X^N)$, 
$$(P^i_tf)(x_1,\dots,x_N) = P^i_t(f(x_1,\dots,x_{i-1},\cdot,x_{i+1},\dots,x_N))(x_i),\qquad 1\le i\le N,$$
we have non-interacting particles.  When the semigroups $P^i_t$, $t\ge0$, $1\le i\le N$, are the same, we have identical non-interacting particles.
If $V(x_1,\dots,x_N)$ is a potential in $C(X^N)$, particle interactivity is then 
modelled by the Schrodinger semigroup $P_t^V$, $t\ge0$, on $C(X^N)$. 

If (A) holds for single particle Markov semigroups $P^i_t$, $t\ge0$, $1\le i\le N$, on $C(X)$, then (A) holds
(with $\epsilon$ replaced by $\epsilon^N$) for the product Markov semigroup $P_t$, $t\ge0$, on $C(X^N)$, corresponding
to non-interacting particles. Similarly for (B). If (C) and (D) hold for $P^i_t$, $t\ge0$, $1\le i\le N$, on $C(X)$, then (C) and (D) hold for the product Markov semigroup $P_t$, $t\ge0$, on $C(X^N)$, assuming $\mathcal D^\infty(X^N)$ can be chosen to be a tensor product of $\mathcal D^\infty(X)$ in a suitable sense. This is the case for the examples above when $\mathcal D^\infty(X^N)=C^\infty(X^N)$ and $\mathcal D^\infty(X)=C^\infty(X)$.

A potential $V$ in $C(X^N)$ is separable  if it is of the form (\ref{eq:ext}) for some $v$ in $C(X)$. 
We are interested in Schrodinger semigroups on $C(X^N)$ with generators of the form $L+V_0+V$ with $V_0,V$ in $C(X^N)$ and $V$ separable. 

Given $f\in C(X^N)$ and a permutation $\sigma$ of $(1,\dots,N)$, let 
$$f^\sigma(x_1,\dots,x_N)=f(x_{\sigma1},\dots,x_{\sigma N}).$$
Given a measure $\mu$ on $X^N$, let $\mu^\sigma$ be the measure with action $\mu^\sigma(f)=\mu(f^\sigma)$. 
A potential $V$ on $X^N$ is symmetric if $V^\sigma=V$ and a measure $\mu$ on $X^N$ is symmetric if $\mu^\sigma=\mu$, both for all permutations $\sigma$. 

Let $P_t$, $t\ge0$, be a Markov semigroup on $C(X^N)$ with generator $L$. We say the semigroup $P_t$, $t\ge0$, is symmetric if $(P_tf)^\sigma=P_t f^\sigma$, $t\ge0$, for all permutations $\sigma$. When the semigroup is symmetric and $V$ is symmetric, we can restrict the supremum in (\ref{eq:DV}) (with $X$ replaced by $X^N$) to symmetric measures. Note for $\mu$ symmetric with marginal $\rho$ and $V$ separable, we have $\mu(V)=\rho(v)$.

Here is the Hohenberg-Kohn theorem in this setting.

\begin{theorem}
\label{thm:HK}
Let $P_t$, $t\ge0$ be a Markov semigroup  on $C(X^N)$  satisfying (A), (B), (C), (D) and let $V_0$ be a 
potential and $V_1$, $V_2$ separable potentials, all in $C(X^N)$, with $V_1$, $V_2$, arising from $v_1$, $v_2$ in 
$C(X)$. Assume $V_0+V_1$ and $V_0+V_2$ are smooth. Let $\mu_1$, $\mu_2$ be symmetric equilibrium 
measures for $V_0+V_1$, $V_0+V_2$ and let $\rho_1$, $\rho_2$ denote their $1$-particle marginals.
Then $\rho_1=\rho_2$ implies $v_1-v_2$ is constant.
\end{theorem}

For example this applies if $V_0$ is symmetric and $P_t$, $t\ge0$, corresponds to non-interacting identical particles.

The proof of this is so short we present it right away.

\proof[Proof of Theorem \ref{thm:HK}]

If $\mu_1$ is an equilibrium measure for $V_0+V_2$, then by Theorem \ref{thm:N=1}, $V_1-V_2 = (V_0+V_1) - (V_0+V_2)$ 
is constant on $X^N$, but $V_1-V_2$ is separable, hence $v_1-v_2$ is constant on $X$. Otherwise, we have
$$\mu_1(V_0+V_2)-I(\mu_1) < \lambda_{V_0+V_2} = \lambda_{V_0+V_2} - \lambda_{V_0+V_1} + \mu_1(V_0+V_1)-I(\mu_1)$$
which implies
$$\rho_1(v_2-v_1) = \mu_1(V_2-V_1) <  \lambda_{V_0+V_2} - \lambda_{V_0+V_1}$$
hence
$$\rho_1(v_2-v_1) <  \lambda_{V_0+V_2} - \lambda_{V_0+V_1}.$$
Reversing the roles of $V_1$, $V_2$,
$$\rho_2(v_1-v_2) <  \lambda_{V_0+V_1} - \lambda_{V_0+V_2}.$$
Since $\rho_1=\rho_2$, this is a contradiction.
\qed

Let $I(\mu)$ correspond to a symmetric Markov semigroup on $C(X^N)$, and let $V_0$, $V$ be in $C(X^N)$ with $V_0$ symmetric and $V$  separable. Let
$$I_{HK}(\rho) \equiv \inf_{\mu\to\rho}\left(I(\mu)-\int_{X^N}V_0\,d\mu\right),$$
where the infimum is over all symmetric $\mu$ in $M(X^N)$ with marginal $\rho$ in  $M(X)$.
Then (\ref{eq:DV}) written over $M(X^N)$ reduces to
$$\lambda_{V_0+V} =
\sup_\mu\left(\int_{X^N}(V_0+V)\,d\mu - I(\mu)\right) =
\sup_\rho \left(\int_Xv\,d\rho - I_{HK}(\rho)\right).$$
Thus the computation of the principal eigenvalue is reduced to computing the $M(X^N)$ universal object 
$I_{HK}$ followed by an optimization over $M(X)$. In the electronic case, density functional theory is the study of approximations of $I_{HK}$ \cite{L}, \cite{P}.

The following sections contain the proofs of Theorems \ref{theorem:equiv}, \ref{theorem:pf}, \ref{thm:N=1}
and supporting Lemmas. Many of the Lemmas are basic and go back to the early papers \cite{DV}, \cite{DV1} and the book \cite{DS}.

\section{The Schrodinger semigroup}

Let $X$ be a compact metric space, let $C(X)$ denote the space of real continuous functions with the sup
norm $\norm{\cdot}$, and let $M(X)$ denote the space of Borel probability measures with the topology of weak convergence. Then $M(X)$ is a compact metric space. Throughout $\mu(f)$ denotes the integral of $f$ against $\mu$.

A strongly continuous positive semigroup on $C(X)$ is a semigroup $P_t$, $t\ge0$,  of bounded operators on $C(X)$ preserving positivity $P_tf\ge0$, for $f\ge0$, $t\ge0$,  and satisfying $\norm{P_tf-f}\to0$ as $t\to0+$.
A Markov semigroup on $C(X)$ is a strongly continuous positive semigroup on $C(X)$ satisfying $P_t1=1$, $t\ge0$. 

Let $C^+(X)$ the strictly positive functions in $C(X)$. Then  $P_tf\in\ C^+(X)$ when $f\in C^+(X)$. 

The subspace $\mathcal D\subset C(X)$ of functions $f\in C(X)$ for which the limit
\begin{equation}
\label{eq:gen}
\lim_{t\to0+}\frac1t \left(P_tf-f\right)
\end{equation}
exists in $C(X)$ is dense. If $Lf$ is defined to be this limit, then $P_t(\mathcal D)\subset  \mathcal D$, $t\ge0$,  the $C(X)$-valued map $t\mapsto P_tf$ is differentiable on $(0,\infty)$ for $f\in\mathcal D$,  and $(d/dt)P_tf=L(P_tf)=P_t(Lf)$, for $f\in\mathcal D$ and $t>0$.

Given $V$ in $C(X)$, the Schrodinger semigroup may be constructed as the unique solution $u(t)=P_t^Vf$, $t\ge0$, of 
\beq{eq:schro}
u(t)=P_tf+\int_0^t P_{t-s}Vu(s)\,ds,\qquad t\ge0.
\eeq
for $f\in C(X)$. 
Then $P^V_t$, $t\ge0$, is a strongly continuous positive semigroup on $C(X)$, and the limit
\beq{eq:pv}
\lim_{t\to0+}\frac1t \left(P^V_tf-f\right)
\eeq
exists in $C(X)$ if and only if $f\in\mathcal D$, in which case it equals $(L+V)f$.  Moreover $P^V_t(\mathcal D)\subset  \mathcal D$, $t\ge0$,  the $C(X)$-valued map $t\mapsto P^V_tf$ is differentiable on $(0,\infty)$ for $f\in\mathcal D$, and  $(d/dt)P^V_tf = (L+V)(P^V_tf) = P^V_t(Lf+Vf)$, for $f\in\mathcal D$ and $t>0$.

For $f\ge0$, (\ref{eq:schro}) implies
\beq{eq:major}
e^{t\min V} P_tf\le P^V_tf\le e^{t\max V}P_tf,\qquad t\ge0.
\eeq
This implies  
$$\min V\le \lambda_V\le \max V.$$

Let $\mathcal D^+$ be the strictly positive functions in $\mathcal D$. For $\mu$ in $M(X)$, let
$$I^V(\mu)\equiv I(\mu)-\int_XV\,d\mu+\lambda_V=-\inf_{u\in\mathcal D^+}\int_X\frac{(L+V-\lambda_V)u}{u}\,d\mu$$
Then $I^0(\mu)=I(\mu)$ and  $I^V(\mu)=0$ iff $\mu$ is an equilibrium measure for $V$.

\begin{lemma}
\label{lemma:lsc}
For $V$ in $C(X)$, $I^V$ is lower semicontinuous, convex, and $0\le I^V \le +\infty$. 
In particular, $I$ is lower semicontinuous, convex, and $0\le I \le +\infty$.
\end{lemma}

\proof
Lower semicontinuity and convexity follow from the fact that $I^V$ is the supremum of continuous affine functions.
The Donsker-Varadhan formula implies $I^V$ is nonnegative.\qed

\begin{lemma}
\label{lemma:core}
Let $\mathcal D^\infty$ be a core for $P_t$, $t\ge0$, that is closed under multiplication. If $f,g\in \mathcal D^\infty$, (\ref{eq:gammaineq}) holds.
\end{lemma}

\proof Expanding both sides of
$$\int_Xp(t,x,dy)f(y)\left(g(y)-g(x)\right)^2 \ge \min f \cdot  \int_Xp(t,x,dy)\left(g(y)-g(x)\right)^2$$
 yields
$$P_t(fg^2)-2gP_t(fg)+g^2P_tf \ge \min f \cdot \left(P_t(g^2) - 2gP_tg+g^2\right)$$
hence
\beqan
(P_t(fg^2)-fg^2)&-&2g(P_t(fg)-fg)+g^2(P_tf-f) \nonumber \\
&\ge& \min f \cdot \left( (P_t(g^2)-g^2) -2g(P_tg-g)\right). \nonumber
\eeqan
Dividing by $t$ and sending $t\to0$ yields half the result. The other half is obtained by replacing $f$ by $-f$.
\qed

Note when $P_t$, $t\ge0$, is a diffusion, e.g. our first example above, one has $L(fg^2)-2gL(fg)+g^2Lf = f\cdot \Gamma(g)$.

For $t>0$ and $u$ in $C^+(X)$, (\ref{eq:major}) implies
$$\log\left(\frac{e^{-\lambda_Vt}P^{V}_tu}{u}\right)$$
is in $C(X)$. 

\begin{lemma}For $V$ in $C(X)$, $\mu\in M(X)$, and $u$ in $C^+(X)$,
\label{lemma:Iapprox}
\beq{eq:log}
\int_X\log\left(\frac{e^{-\lambda_Vt}P^{V}_tu}{u}\right)\,d\mu\ge-tI^V(\mu),\qquad t\ge0.
\eeq
\end{lemma}

The proof follows that of Lemma 3.1 in \cite{DV1}.

\proof
By definition of $I^V(\mu)$,
\beq{eq:equil}
\int_X\frac{(L+V-\lambda_V)u}{u}\,d\mu\ge-I^V(\mu),\qquad u\in\mathcal D^+.
\eeq
When $I^V(\mu)=+\infty$, the result is valid, hence we may assume $I^V(\mu)<\infty$. 
For $t=0$, (\ref{eq:log}) is an equality. Moreover for $t>0$ and $u\in\mathcal D^+$, by (\ref{eq:major}) we have $e^{-\lambda_Vt}P^V_tu\in\mathcal D^+$ and
$$\frac d{dt} \int_X\log\left(\frac{e^{-\lambda_Vt}P^{V}_tu}{u}\right)\,d\mu=\int_X\frac{(L+V-\lambda_V)(e^{-\lambda_Vt}P^V_tu)}{e^{-\lambda_Vt}P^V_tu}\,d\mu\ge -I^V(\mu).$$
This establishes (\ref{eq:log}) for $u\in\mathcal D^+$. Since $\mathcal D^+$ is dense in $C^+(X)$, (\ref{eq:log}) is valid for $u$ in $C^+(X)$. \qed

\section{Equilibrium Measures}

Let $L^1(\mu)$ denote the $\mu$-integrable Borel functions on $X$ with  
$$\norm{f}_{L^1(\mu)}=\int_X|f|\,d\mu=\mu(|f|).$$

The following strengthening of Lemma \ref{lemma:Iapprox}  is necessary in the next section. Let $B(X)$ denote the bounded Borel functions on $X$. Recall (\ref{eq:posmeas}) $0\le P^V_tu(x)\le +\infty$ is well-defined for $u\ge0$ Borel, for all $x\in X$.

\begin{lemma} 
\label{lemma:strong}
Fix $V\in C(X)$ and $\mu\in M(X)$.  Let $u>0$ Borel satisfy $\log u\in L^1(\mu)$. Then for $t\ge0$,
\beq{eq:logP}
tI^V(\mu)+\int_X\log^+\left(\frac{e^{-\lambda_Vt}P^{V}_tu}{u}\right)\,d\mu\ge\int_X\log^-\left(\frac{e^{-\lambda_Vt}P^{V}_tu}{u}\right)\,d\mu.
\eeq
Here the integrals may be infinite.
\end{lemma}

\proof We may assume $I^V(\mu)<\infty$, otherwise  (\ref{eq:logP}) is true.

Let $u>0$ be Borel with $\log u \in L^1(\mu)$. We establish (\ref{eq:logP}) in three stages, first for $\log u\in B(X)$, then for $\log u$ bounded below, then in general.
Let $Q_t=e^{-\lambda_Vt}P^V_t$, $t\ge0$. 

Suppose $|\log u| \le M$ and suppose $u_n>0$, $n\ge1$,  satisfy $|\log u_n| \le M$, $n\ge1$. If $u_n\to u$ pointwise on $X$, it follows that $Q_tu_n\to Q_t u$ pointwise on $X$. Assume  (\ref{eq:log}) is valid for $u_n$, $n\ge1$.
Since by (\ref{eq:major})
$$t(\min V-\lambda_V) -2M \le \log \left(\frac{Q_t u_n}{u_n}\right)  \le t(\max V-\lambda_V) + 2M,\qquad n\ge1,$$
it follows that (\ref{eq:log}) is valid for $u$. Thus the set of Borel $f$ in $B(X)$ with $u=e^f$ satisfying (\ref{eq:log}) is closed under bounded pointwise convergence. Since (\ref{eq:log}) is valid when $f=\log u\in C(X)$, it follows that (\ref{eq:log}) hence (\ref{eq:logP}) is valid for all Borel $u$ satisfying $\log u\in B(X)$.
 Here both sides of (\ref{eq:logP}) are finite.

Next, assume $\log u$ in $L^1(\mu)$ and $u\ge\delta>0$ and let $u_n=u\wedge n$, $n\ge1$. Then
$$\log\left(\frac{Q_tu}{u}\right) \ge \log\left(\frac{Q_tu_n}{u}\right)
= \log\left(\frac{Q_tu_n}{u_n}\right) + \log\left(\frac{u_n}{u}\right)$$
so
$$\log^+\left(\frac{Q_tu}{u}\right) \ge \log^-\left(\frac{Q_tu}{u}\right) +
\log\left(\frac{Q_tu_n}{u_n}\right) + \log\left(\frac{u_n}{u}\right).$$
Hence
$$\int_X \log^+\left(\frac{Q_tu}{u}\right)\,d\mu \ge
\int_X \log^-\left(\frac{Q_tu}{u}\right)\,d\mu - tI^V(\mu) + \int_{u>n}(\log n - \log u)\,d\mu.$$
Discarding the $\log n$ term and passing to the limit $n\to\infty$ yields (\ref{eq:logP}). Note $u\ge\delta$ and (\ref{eq:major}) imply 
$$\log^-\left(\frac{Q_tu}{u}\right)
= \log^+\left(\frac{u}{Q_tu}\right)
\le |\log u| + (\lambda_V-\min V)t + \log\frac1\delta$$
so the right side of (\ref{eq:logP}) is finite in this case and in fact (\ref{eq:log}) is valid.

Now assume $\log u$ in $L^1(\mu)$ and let $u_\delta=u\vee\delta$. Then 
$$\log^+\left(\frac{Q_tu_\delta}{u}\right) = \log^-\left(\frac{Q_tu_\delta}{u}\right) +
\log\left(\frac{Q_tu_\delta}{u_\delta}\right) + \log\left(\frac{u_\delta}{u}\right)$$
so
$$\int_X \log^+\left(\frac{Q_tu_\delta}{u}\right)\,d\mu \ge
\int_X \log^-\left(\frac{Q_tu_\delta}{u}\right)\,d\mu - tI^V(\mu) + 
\int_{u<\delta}\log\left(\frac{\delta}{u}\right)\,d\mu$$
hence
\beq{eq:logPdelta}
tI^V(\mu) + \int_X\log^+\left(\frac{Q_tu_\delta }{u}\right)\,d\mu  
\ge \int_X\log^-\left(\frac{Q_tu_\delta}{u}\right)\,d\mu,
\eeq
where we discarded the right-most integral as its integrand is nonnegative.
To establish (\ref{eq:logP}), we pass to the limit $\delta\downarrow0$ in (\ref{eq:logPdelta}). 
We may assume
$$\int_X\log^+\left(\frac{Q_tu}{u}\right)\,d\mu<\infty,$$
otherwise (\ref{eq:logP}) is true. This implies $\log^+(Q_tu/u)(x)<\infty$ for $\mu$-a.a $x$ which implies $Q_tu(x)<\infty$ for $\mu$-a.a. $x$.  Since $u_\delta\le u+1$ for $\delta<1$, it follows by the dominated convergence theorem that $Q_tu_\delta\to Q_tu$ a.e. $\mu$ as $\delta\downarrow0$. 

Since
$$\log^-\left(\frac{Q_tu_\delta}{u}\right),\qquad \delta>0,$$
increases as $\delta\downarrow0$, the right side of (\ref{eq:logPdelta}) converges to the right side of (\ref{eq:logP}). Using $2\log^+(a+b)\le 2\log2 + \log^+a+\log^+b$, (\ref{eq:major}), and $u_\delta\le u+1$ for $\delta<1$, we have
$$2\log^+\left(\frac{Q_tu_\delta }{u}\right)\le 
2\log2 + \log^+\left(\frac{Q_tu }{u}\right) + |\log u| + t(\max V-\lambda_V),$$
hence 
the dominated convergence theorem shows the left side of (\ref{eq:logPdelta}) converges to the left side of (\ref{eq:logP}).
\qed

Let $P^{V,\psi}_t$ be as in (\ref{eq:Ppsi}).

\begin{corollary}
\label{corollary:ineqQ}
Fix $V\in C(X)$, $\mu\in M(X)$, let $\log\psi\in L^1(\mu)$, and let $u>0$ Borel satisfy $\log u\in L^1(\mu)$. Then for $t\ge0$,
\beq{eq:logQ}
tI^V(\mu)+\int_X\log^+\left(\frac{P^{V,\psi}_tu}{u}\right)\,d\mu\ge\int_X\log^-\left(\frac{P^{V,\psi}_tu}{u}\right)\,d\mu.
\eeq
Here the integrals may be infinite.
\end{corollary}

\proof Since $\log\psi$ is in  $L^1(\mu)$, $\log(u\psi)$ is in $L^1(\mu)$ iff $\log u$ is in $L^1(\mu)$. Now apply Lemma
\ref{lemma:strong}.\qed

\begin{corollary}
\label{corollary:nonneg}
Let $V\in C(X)$ and $\log\psi\in L^1(\mu)$. Then $\mu\in M(X)$ is an equilibrium measure for $V$ iff 
 $$
\int_X\log^+\left(\frac{P^{V,\psi}_tu}{u}\right)\,d\mu\ge\int_X\log^-\left(\frac{P^{V,\psi}_tu}{u}\right)\,d\mu
$$
for $t\ge0$ and $u>0$ satisfying $\log u\in L^1(\mu)$.
\end{corollary}

\proof
If $\mu$ is an equilibrium measure, $I^V(\mu)=0$ so the result follows from Corollary \ref{corollary:ineqQ}.
Conversely, assume the inequality holds for all $u>0$ satisfying $\log u\in L^1(\mu)$.
For $u\in C^+(X)$, the function $u/\psi$  satisfies $\log (u/\psi)\in L^1(\mu)$. Inserting $u/\psi$ in the inequality yields
 $$
 \int_X\log^+\left(\frac{e^{-\lambda_Vt}P^V_tu}{u}\right)\,d\mu\ge\int_X\log^-\left(\frac{e^{-\lambda_Vt}P^V_tu}{u}\right)\,d\mu.
$$
For $u$ in $C^+(X)$, the integrals are finite hence
$$
 \int_X\log\left(\frac{e^{-\lambda_Vt}P^V_tu}{u}\right)\,d\mu\ge0.
 $$
For $u\in\mathcal D^+$, with $Q_t=e^{-\lambda_Vt}P^V_t$, $t\ge0$, we have $Q_tu\in\mathcal D^+$ so
$$Q_tu = u + t(L+V-\lambda_V)u + o(t),\qquad t\to0,$$

$$\frac{Q_tu}{u} = 1 + t\frac{(L+V-\lambda_V)u}{u} + o(t),\qquad t\to0,$$

$$\log\left(\frac{Q_tu}{u}\right) =  t\frac{(L+V-\lambda_V)u}{u} + o(t),\qquad t\to0,$$
 all uniformly on $X$. Hence  dividing by $t$ and sending $t\to0$ yields
$$\int_X\frac{(L+V-\lambda_V)u}{u}\,d\mu\ge0. 
$$
This implies $I^V(\mu)\le0$, hence $I^V(\mu)=0$.\qed

A strongly continuous positive semigroup on $L^1(\mu)$ is a semigroup $P_t$, $t\ge0$,  of bounded operators on $L^1(\mu)$ preserving positivity $P_tf\ge0$ a.e. $\mu$, for $f\ge0$ a.e. $\mu$, $t\ge0$,  and satisfying $\norm{P_tf-f}_{L^1(\mu)}\to0$ as $t\to0+$.
A Markov semigroup on $L^1(\mu)$ is a strongly continuous positive semigroup on $L^1(\mu)$ satisfying $P_t1=1$ a.e. $\mu$, $t\ge0$.

\begin{lemma}
\label{lemma:Q}
Let $V\in C(X)$ and suppose $\pi$ and $\mu$ are measures with $\mu<<\pi$, and let $\psi=d\mu/d\pi$. If $\pi$ is a ground
measure for $V$, then $P^{V,\psi}_t|f|(x)<\infty$ for $\mu$-a.a. $x$ and $f$ in $L^1(\mu)$, $P^{V,\psi}_t$, $t\ge0$, is a strongly continuous positive semigroup on $L^1(\mu)$, and
\beq{eq:Qinv}
\mu(P^{V,\psi}_tf)=\mu(f),\qquad t\ge0,
\eeq
for $f$ in $L^1(\mu)$.  If $\psi$ is a ground state for $V$ relative to $\mu$,  $P^{V,\psi}_t$, $t\ge0$, is a Markov semigroup on $L^1(\mu)$.
\end{lemma}

\proof If $\pi$ is a ground measure, for $f$ in $C(X)$ we have
\beqan
\norm{e^{-\lambda_Vt}P^V_tf}_{L^1(\pi)}&=&\int_X|e^{-\lambda_Vt}P^V_tf|\,d\pi \nonumber \\
&\le& \int_Xe^{-\lambda_Vt}P^V_t|f|\,d\pi =\int_X|f|\,d\pi=\norm{f}_{L^1(\pi)}.\nonumber
\eeqan
Hence $e^{-\lambda_Vt}P^V_t$, $t\ge0$, satisfies
\beq{eq:piest}
\norm{e^{-\lambda_Vt}P^V_tf}_{L^1(\pi)} \le \norm{f}_{L^1(\pi)},\qquad t\ge0,
\eeq
for $f$ in $C(X)$. Since the collection of functions $f$ satisfying (\ref{eq:piest}) is closed under bounded pointwise convergence, (\ref{eq:piest}) is valid for $f\in B(X)$. Inserting $f\wedge n$ with $f$ nonnegative Borel and sending $n\to\infty$, (\ref{eq:piest}) is then valid for nonnegative Borel $f$. It follows that $e^{-\lambda_Vt}P^{V}_t|f|(x)<\infty$, $\pi$-a.a. $x$, for $f$ in $L^1(\pi)$, hence $e^{-\lambda_Vt}P^V_t$, $t\ge0$, are well-defined contractions on $L^1(\pi)$.  By (\ref{eq:piest}) and the density of $C(X)$ in $L^1(\pi)$, this implies 
$\pi(e^{-\lambda_Vt}P^V_tf)=\pi(f)$, $t\ge0$, for $f$ in $L^1(\pi)$ and implies $e^{-\lambda_Vt}P^V_t$, $t\ge0$, is a strongly continuous positive semigroup on $L^1(\pi)$. 

Since $\psi\in L^1(\pi)$, (\ref{eq:Qinv}) follows for $f\in C(X)$. But  (\ref{eq:piest}) for $f$ nonnegative Borel implies
\beq{eq:qest}
\norm{P^{V,\psi}_tf}_{L^1(\mu)} \le \norm{f}_{L^1(\mu)},\qquad t\ge0,
\eeq
for $f$ nonnegative Borel, hence $P^{V,\psi}_t|f|(x)<\infty$, $\mu$-a.a. $x$, for $f$ in $L^1(\mu)$, hence $P^{V,\psi}_t$, $t\ge0$, are well-defined contractions on $L^1(\mu)$. Moreover
$$\norm{P^{V,\psi}_tf-f}_{L^1(\mu)}=\norm{e^{-\lambda_Vt}P^V_t(f\psi)-f\psi}_{L^1(\pi)}\to0,\qquad t\to0+,f\in C(X).$$ 
By (\ref{eq:qest}) and the density of $C(X)$ in $L^1(\mu)$, we conclude $P^{V,\psi}_t$, $t\ge0$, is a strongly continuous positive semigroup on $L^1(\mu)$ and (\ref{eq:Qinv})  holds for $f\in L^1(\mu)$. 

If  $\psi$ is a ground state relative to $\mu$, $P^{V,\psi}_t1=1$ a.e. $\mu$. Thus in this case $P^{V,\psi}_t$, $t\ge0$, is a Markov semigroup on $L^1(\mu)$.\qed

\section{Proofs of the Theorems}

\proof[Proof of Theorem \ref{theorem:equiv}] 
For the first assertion, we have a ground measure $\pi$ for $V$  and a ground state $\psi$ for $V$ relative to $\mu$ satisfying $\log\psi\in L^1(\mu)$.
Suppose  $\log u\in L^1(\mu)$. Then $P^{V,\psi}_t|\log u|$ is in $L^1(\mu)$ and there is a set $N$ with $\mu(N)=0$ and $P^{V,\psi}_t(|\log u|)(x)<\infty$ and 
$P^{V,\psi}_t1(x)=1$ for $x\not\in N$. Jensen's inequality applied to the integral 
$f\mapsto (P^{V,\psi}_tf)(x)$ (see (\ref{eq:posmeas})) implies 
$$\log\left(\frac{P^{V,\psi}_tu}{u}\right)(x)\ge P^{V,\psi}_t(\log u)(x)-(\log u)(x),\qquad x\not\in N,$$
hence for $x\not\in N$,
$$
\log^+\left(\frac{P^{V,\psi}_tu}{u}\right)(x)\ge \log^-\left(\frac{P^{V,\psi}_tu}{u}\right)(x) 
+ P^{V,\psi}_t(\log u)(x)-(\log u)(x).
$$
Integrating over $X$ against $\mu$, the integrals of the right-most two terms cancel by (\ref{eq:Qinv}) hence by  Corollary \ref{corollary:nonneg},  $\mu$ is an equilibrium measure for $V$, establishing the first assertion.

For the second assertion, assume $\pi$ is a ground measure for $V$ and $\mu$ is an equilibrium measure for $V$.
Note $\int P^{V,\psi}_t1\,d\mu<\infty$ so $\int \log^+\left(P^{V,\psi}_t1\right)\,d\mu<\infty$. 
By Corollary \ref{corollary:nonneg}, it follows that $\int \log^-\left(P^{V,\psi}_t1\right)\,d\mu<\infty$, hence $\log\left(P^{V,\psi}_t1\right)$ is in $L^1(\mu)$. By Jensen's inequality, (\ref{eq:Qinv}), and Corollary \ref{corollary:nonneg},
$$0=\log(\mu(1))=\log\left(\int_X P^{V,\psi}_t1\,d\mu\right)\ge \int_X \log(P^{V,\psi}_t1)\,d\mu\ge0.$$
Since $\log$ is strictly concave, this can only happen if $P^{V,\psi}_t1$ is $\mu$ a.e. constant. By (\ref{eq:Qinv}), the constant
is $1$. Since $\psi>0$ a.e. $\mu$ is immediate, this establishes the second assertion.

For the third assertion, 
assume $\mu$ is an equilibrium measure for $V$ and $\psi$ is a ground state for $V$ relative to $\mu$.
Then $P^{V,\psi}_t1=1$ a.e. $\mu$, so for $u\in C^+(X)$,
$$\frac{\min u}{\max u}\le \frac{P^{V,\psi}_tu}{u}\le\frac{\max u}{\min u},\qquad a.e. \mu,$$
hence $\log(P^{V,\psi}_tu/u)$ 
  is in $L^1(\mu)$ for $u\in C^+(X)$. By Corollary \ref{corollary:nonneg}, for $f\in C(X)$,
$$\beta(\epsilon)\equiv\int_X\log\left(\frac{P^{V,\psi}_te^{\epsilon f}}{e^{\epsilon f}}\right)\,d\mu\ge0,\qquad |\epsilon|<1,$$
and $\beta(0)=0$, hence  $\dot\beta(0)=0$.
Differentiating at $\epsilon=0$, we obtain
\beq{inv1}
\int_X e^{-\lambda_Vt}P^V_t(f\psi)\,d\pi = \int_X f\psi\, d\pi
\eeq
for $f\in C(X)$. Since the collection of functions $f$ satisfying (\ref{inv1}) is closed under bounded pointwise convergence,  (\ref{inv1}) holds for $f\in B(X)$. Now for $f\in C(X)$, $f_\epsilon \equiv f\psi/(\psi+\epsilon)\to f$ boundedly as $\epsilon\downarrow0$, thus replacing $f$ by $f/(\psi+\epsilon)$ in (\ref{inv1}) and letting $\epsilon\downarrow0$ 
establishes (\ref{eq:ground}), hence $\pi$ is a ground measure for $V$.
This establishes the third assertion.\qed

For $\mu,\pi$ in $M(X)$, the entropy  of $\mu$ relative to $\pi$ is
$$H(\mu,\pi)\equiv\sup_V\left(\int_XV\,d\mu-\log\int_Xe^V\,d\pi\right)$$
where the supremum is over $V$ in $C(X)$. 

\begin{lemma} $H(\mu,\pi)\ge0$ is finite iff $\mu<<\pi$ and $\psi=d\mu/d\pi$ satisfies  $\log\psi \in L^1(\mu)$,  in which case
$$H(\mu,\pi)=\int_X\log\psi\,d\mu = \int_X \psi\log \psi\,d\pi.$$
Moreover $H$ is lower-semicontinuous and convex separately in each of $\mu$ and $\pi$.
\end{lemma}

This is Lemma 2.1 in \cite{DV1}.
\proof
The lower-semicontinuity and convexity follow from the definition of $H$ as a supremum of convex functions, in each variable 
$\pi$, $\mu$ separately. Suppose $H(\mu,\pi)<\infty$. Since the set of $V$ in $B(X)$ satisfying
$$\int_XV\,d\mu-\log\int_Xe^V\,d\pi\le H(\mu,\pi)$$
 contains $C(X)$ and is closed under bounded pointwise convergence, it equals $B(X)$. 
Insert $V=r1_A$ into the definition of $H$, where $\pi(A)=0$, obtaining
$$r\mu(A)\le r\mu(A)-\log(\pi(A^c))\le H(\mu,\pi).$$
Let $r\to\infty$ to conclude $\mu<<\pi$. Since $\psi=d\mu/d\pi\in L^1(\pi)$, let $0\le f_n\in C(X)$ with $f_n\to \psi$
in $L^1(\pi)$. By passing to a subsequence, assume $f_n\to \psi$ a.e. $\pi$.
Insert $V=\log(f_n+\epsilon)$ into the definition of $H$ to yield
$$\int_X\log(f_n+\epsilon)\,d\mu-\log\int_X(f_n+\epsilon)\,d\pi\le H(\mu,\pi).$$
Let $n\to\infty$; by Fatou's lemma,
$$\int_X \psi\log(\psi+\epsilon)\,d\pi-\log\int_X(\psi+\epsilon)\,d\pi
\le H(\mu,\pi).$$
Since $\pi(\psi+\epsilon)=1+\epsilon$, applying Fatou's lemma again as $\epsilon\to0$, $\int_X \psi\log \psi\,d\pi\le H(\mu,\pi)$.

Conversely, suppose $\psi=d\mu/d\pi$ exists and $\psi\log \psi\in L^1(\pi)$. By Jensen's inequality,
$$\int_X V\,d\mu\le \log\int_X e^V\,d\mu,\qquad V\in B(X).$$
Replace $V$ by $V-\log(\psi\wedge n+\epsilon)$ to get
$$\int_XV\,d\mu-\log\int_X\left(\frac{e^V\psi}{\psi\wedge n+\epsilon}\right)\,d\pi\le \int_X \psi\log (\psi\wedge n+\epsilon)\,d\pi.$$
Let $\epsilon\to0$ followed by $n\to\infty$ obtaining
$$\int_XV\,d\mu-\log\int_X e^V\,d\pi\le \int_X \psi\log \psi\,d\pi.$$
Now maximize over $V$ in $C(X)$ to conclude $H(\mu,\pi)\le \int_X \psi\log \psi\,d\pi$.
\qed

\proof[Proof of Theorem \ref{theorem:pf}]

By (\ref{eq:log}),
$$\int_X\log\left(\frac{e^{-\lambda_Vt}P^{V}_tu}{u}\right)\,d\mu\ge-tI^V(\mu),\qquad u\in C^+(X).$$
Thus for $f\in C(X)$,
$$\int_Xf\,d\mu-\int_X\log\left(e^{-\lambda_Vt}P^V_te^f\right)\,d\mu\le tI^V(\mu),\qquad f\in C(X).$$
By Jensen's inequality,
$$\int_Xf\,d\mu-\log\int_X\left(e^{-\lambda_Vt}P^V_te^f\right)\,d\mu\le tI^V(\mu),\qquad f\in C(X).$$
Defining
$$\mu_t(f)=e^{-\lambda_Vt}\mu(P^V_tf)$$
and 
$$\pi_t(f)=\frac{\mu_t(f)}{\mu_t(1)}$$
yields
$$\int_Xf\,d\mu-\log\int_Xe^f\,d\pi_t\le tI^V(\mu)+\log \mu_t(1),\qquad f\in C(X).$$
Taking the supremum over all $f$ yields
$$H\left(\mu,\pi_t\right)\le tI^V(\mu)+\log \mu_t(1).$$
Note $\mu_t(1)\le C$, $t\ge0$, hence
$$H\left(\mu,\pi_t\right)\le tI^V(\mu)+\log C,\qquad t\ge0.$$
Now set
$$\bar\pi_T=\frac{\int_0^T \mu_t\,dt}{\int_0^T \mu_t(1)\,dt}=\frac{\int_0^T \mu_t(1)\pi_t\,dt}{\int_0^T\mu_t(1)\,dt},\qquad T>0.$$
Then $\pi_t$ is in $M(X)$ for $t>0$, $\bar\pi_T$ is in $M(X)$ for $T>0$. 

Now assume $\mu$ is an equilibrium measure for $V$; then $I^V(\mu)=0$. By convexity of $H$.
$$H\left(\mu,\bar\pi_T\right)\le \log C,\qquad T>0.$$
By compactness of $M(X)$, select a sequence $T_n\to\infty$ with $\pi_n=\bar\pi_{T_n}$ converging to some $\pi$. By lower-semicontinuity of $H$, we have $H(\mu,\pi)\le \log C$. Thus $\mu<<\pi$ with $\psi=d\mu/d\pi$ satisfying $\psi\log\psi\in L^1(\pi)$. Since
$$\log\mu(e^{-\lambda_Vt}P^V_t1)\ge \mu(\log(e^{-\lambda_Vt}P^V_t1))\ge0,$$
we have $\mu_t(1)\ge1$, $t\ge0$. This is enough to show
$$\pi_n\left(e^{-\lambda_VT}P^V_Tf\right)=\pi_n(f)+o(1),\qquad n\to\infty,$$
for all $T>0$. Thus $\pi$ is a ground measure for $V$. By Theorem \ref{theorem:equiv}, $\psi$ is a ground state for $V$ relative to $\mu$. The
remaining assertions are in Lemma \ref{lemma:Q}.
\qed

We establish two lemmas used in the proof of Theorem \ref{thm:N=1}.

\begin{lemma}
Let $V\in C(X)$. Under assumption (A), (\ref{eq:supbound}) holds.
\end{lemma}

This is Lemma 4.3.1 in \cite{DS}.

\proof
Let $T>0$ and $\epsilon>0$ be as in (A). By (\ref{eq:major}), for $t\ge0$,
\beqan
P_TP^V_t1 &\le& e^{-T\min V} P_T^VP_t^V1 = e^{-T\min V} P_t^VP_T^V1 \nonumber \\
&\le& e^{T(\max V-\min V)} P_t^VP_T1 = e^{T(\max V-\min V)} P_t^V1. \nonumber
\eeqan
Similarly, one has
$$
P_TP^V_t1 \ge e^{T(\min V-\max V)} P_t^V1
$$
hence
$$e^{T(\max V-\min V)} P_t^V1 \ge P_TP^V_t1 \ge e^{T(\min V-\max V)} P_t^V1.$$
Let $\epsilon' =\epsilon e^{2T(\min V-\max V)}$. By (A) this implies
$$P^V_t1(x) \ge \epsilon'  P^V_t1(y),\qquad x,y\in X,$$
 hence
$$\norm{P^V_t} = \sup_x P^V_t1(x) \ge \phi(t)\equiv \inf_x P^V_t1(x) \ge  \epsilon'  \norm{P^V_t},\qquad t\ge0.$$
But $\phi(t)$ is supermultiplicative so
$$\sup_{t>0} \frac1t \log\phi(t) = \lim_{t\to\infty} \frac1t \log \phi(t) \le \lim_{t\to\infty} \frac1t \log\norm{P^V_t} =  \lambda_V.$$
Since $\epsilon'\norm{P^V_t} \le \phi(t)$, this  implies (\ref{eq:supbound})  with $C\le 1/\epsilon'$.\qed

\begin{lemma}
\label{lemma:bounded}
Under assumption (A), the ground state $\psi$ in Theorem \ref{theorem:pf} may be chosen such that $\log\psi$ is in $B(X)$. If moreover (B) holds, $\supp(\mu)=X$. If moreover (C) holds and $V$ is smooth, $\psi$ may be chosen in $\mathcal D^\infty$ and strictly positive, and satisfies
$$L\psi + V\psi = \lambda_V\psi.$$
\end{lemma}

\proof
With $T$ and $\epsilon$ as in (A), let $Q_T=e^{-\lambda_VT}P^V_T$ and $\epsilon' = \epsilon e^{T(\min V-\max V)}$. Then $Q_T\psi = \psi$ a.e. $\mu$. By (A) and (\ref{eq:major}) we have 
\beq{eq:ineq}
Q_T|f|(x) \ge\epsilon' Q_T|f|(y),\qquad x,y\in X,
\eeq
for all $f\in C(X)$. Since the collection of functions $f$ satisfying (\ref{eq:ineq}) is closed under bounded pointwise convergence, (\ref{eq:ineq}) is valid for $f\in B(X)$. 
Hence
$$Q_T\psi(x) \ge Q_T(\psi\wedge n)(x) \ge \epsilon' Q_T(\psi\wedge n)(y),\qquad x,y\in X.$$
Let $\tilde\psi\equiv Q_T\psi$. Sending $n\to\infty$ yields
\beq{eq:everywhere}
\tilde\psi(x)  \ge \epsilon' \tilde\psi(y),\qquad x,y\in X.
\eeq
Since $\psi$ is a ground state, $\tilde\psi=\psi$ a.e. $\mu$. Since $0<\psi<\infty$ a.e. $\mu$, we have $0<\tilde\psi<\infty$ a.e. $\mu$ hence (\ref{eq:everywhere}) implies $\tilde\psi$ is bounded away from zero and away from infinity, i.e.  $\log\tilde\psi$ is in $B(X)$. Since $d\pi=d\mu/\psi = d\mu/\tilde\psi$, Theorem \ref{theorem:equiv} implies $\tilde\psi$ is a ground state. Thus we may replace $\psi$ by $\tilde\psi$ and assume $\log\psi\in B(X)$.

With $T>0$ as in (B), $f\in C(X)$ nonnegative implies 
$$\mu(f) = \mu(P^{V,\psi}_Tf) \ge \frac{\inf\psi}{\sup\psi}e^{T(\min V-\lambda_{V})}\mu(P_Tf) >0.$$
Hence $\supp(\mu)=X$. 

Now let $\tilde\psi\equiv Q_T\psi$ and assume $V$ is smooth. Then $\tilde\psi = \psi$ a.e. $\mu$ 
hence as before $\tilde\psi$ is a ground state. Since $\tilde\psi\in\mathcal D^\infty$,  we may 
replace $\psi$ by $\tilde\psi$ and assume $\psi\in\mathcal D^\infty$. 

Since $\supp(\mu)=X$,  $e^{-\lambda_Vt}P^V_t\psi = \psi$, $t\ge0$, holds identically on $X$, hence $\psi$ is strictly positive.  Differentiating this yields  $L\psi+V\psi=\lambda_V\psi$.\qed

\proof[Proof of Theorem \ref{thm:N=1}]
Let $\psi_i\in\mathcal D^\infty$ be the strictly positive ground states for $V_i$ relative to $\mu$, $i=1,2$, given by Lemma \ref{lemma:bounded}. Since $\mu$ is $P^{V_i,\psi_i}_t$-invariant, $i=1,2$, differentiating (\ref{eq:Qinv}) yields
$$\int_X\frac{L(\psi_i f)}{\psi_i} + V_if - \lambda_{V_i}f \,d\mu = \int_X \left(\frac{L(\psi_i f)}{\psi_i} -f\frac{L\psi_i}{\psi_i}\right)\,d\mu 
= 0,\qquad f\in\mathcal D^\infty$$
for $i=1,2$. Subtract these two equations and insert $f=\psi_1/\psi_2$ to get
$$\int_X \left(\frac{L(fg^2)-2gL(fg)+g^2Lf}{fg}\right)\,d\mu = 0,$$
where now $f=\psi_2$ and $g=\psi_1/\psi_2$. But by Lemma \ref{lemma:core}, 
$$\frac{L(fg^2)-2gL(fg)+g^2Lf}{fg} \ge \frac{\min f}{\max fg} \Gamma(g),$$
so
$$\int_X\Gamma(\psi_1/\psi_2)\,d\mu =0.$$
Since $\supp(\mu)=X$,  $\Gamma(\psi_1/\psi_2)\equiv0$  which by (D) yields $\psi_1=c\psi_2\equiv\psi$. Thus we arrive at
$L\psi+V_i\psi=\lambda_{V_i}\psi$ for $i=1,2$. Subtracting yields the result.\qed


\begin{bibdiv}
\begin{biblist}

\bib{A}{article}{
	title={Uniform Positivity Improving Property, Sobolev Inequalities, and Spectral Gaps},
	author={S. Aida},
	journal={J. Functional Analysis},
	volume={158},
	date={1998},
	pages={152--185}
}

\bib{B}{incollection}{
	title={Functional inequalities for Markov semigroups},
	author={D. Bakry},
	date={2004},
	booktitle={Probability measures on groups},
}

\bib{DS}{book}{
	title={Large Deviations},
	series={Pure and Applied Mathematics Series},
	publisher={Academic Press},
	author={J.-D. Deuschel},
	author={D. W. Stroock},
	date={1984},
	volume={137}
}

\bib{DV}{article}{
        title={On a variational formula for the principal eigenvalue for operators with maximum principle},
        author={M. D. Donsker},
	author={S. R. S. Varadhan},
        journal={Proceedings of the National Academy of Sciences USA},
        volume={72},
        date={1975},
        pages={780--783},
	url={http://www.pnas.org/content/72/3/780.full.pdf}
}

\bib{DV1}{article}{
	author={M. D. Donsker},
	author={S. R. S. Varadhan},
	date={1975},
	title={Asymptotic evaluation of certain Markov process expectations for large time, I,},
	journal={Communications on Pure and Applied Mathematics}, 
	volume={\bf XXVIII},
	pages={1--47}
}

\bib{F}{article}{
	author={S. Friedland},
	date={1981},
	title={Convex spectral functions},
	journal={Linear and Multilinear Algebra},
	volume={9}, 
	pages={299--316}
}


\bib{G}{article}{
	author={L. Gross},
	date={1972},
	title={Existence and uniqueness of physical ground states},
	journal={J. Func. Anal.},
	volume={10},
	pages={52--109}
}

\bib{HK}{article}{
	author={P. Hohenberg},
	author={W. Kohn},
	date={1964},
	title={Inhomogeneous Electron Gas},
	journal={Phys. Rev. B}, 
	volume={136}, 
	pages={864--871}
}

\bib{L}{article}{
	author={E. H. Lieb},
	date={1983},
	title={Density Functionals for Coulomb Systems},
	journal={International Journal for Quantum Chemistry},
	volume={XXIV}, 
	pages={243--277}
}

\bib{P}{article}{
	author={J.P. Perdew},
	author={K. Burke},
	author={M. Ernzerhof},
	date={1996},
	title={Generalized Gradient Approximation Made Simple},
	journal={Phys. Rev. Lett.},
	volume={77}, 
	pages={3865--3868}
}

\bib{SI}{article}{
	author={B. Simon},
	date={1982},
	title={Schrodinger Semigroups},
	journal={Bulletin AMS}, 
	volume={7}, 
	pages={447-526}
}



\end{biblist}

\end{bibdiv}
\end{document}